\documentclass{amsproc}
\usepackage{amssymb}
\usepackage{amsmath}

\newtheorem{theorem}{Theorem}[section]
\newtheorem{lemma}[theorem]{Lemma}
\newtheorem{proposition}[theorem]{Proposition}
\newtheorem{corollary}[theorem]{Corollary}
\theoremstyle{definition}

\theoremstyle{remark}
\newtheorem{remark}[theorem]{Remark}

\numberwithin{equation}{section}

\newcommand{\norm}[1]{\left\Vert#1\right\Vert}
\newcommand{\set}[1]{\left\{#1\right\}}

\begin{document}

\title{On $C^*$-algebras generated by some deformations of CAR relations}

\author{Daniil Proskurin}
\address{Department of Cybernetics, Kyiv Taras Shevchenko University,
 Volodymyrska, 64, 01033 Kyiv, Ukraine}
\email{prosk@univ.kiev.ua.}
\thanks{The work was completed when the first author was
visiting Chalmers University of Technology in G{\"o}teborg, Sweden.
The visit was supported by a grant from the Swedish Royal Academy of
Sciences as a part of program of cooperation with former Soviet
Union and by STINT.}

\author{Yuri\u\i{} Savchuk}
\address{Department of Cybernetics, Kyiv Taras Shevchenko University,
 Volodymyrska, 64, 01033 Kyiv, Ukraine}
\email{yurius@univ.kiev.ua}
\author{Lyudmila Turowska}
\address{Chalmers Tekniska H{\"o}gskola, Matematiska Vetenskaper, SE 412 96
G{\"o}teborg, Sweden}
\email{turowska@math.chalmers.se}
\thanks{Lyudmila Turowska was partially supported by
Swedish Research Council.}
\date{}



\begin{abstract}
We study the representations and enveloping $C^*$-algebras for Wick
analogues of CAR and twisted CAR algebras. Realizations of the
considered $C^*$-algebras are given as algebras of continuous
matrix-functions satisfying certain boundary conditions.
\end{abstract}

\maketitle

\section*{Introduction}
In this paper we study $*$-representations and enveloping
$C^*$-algebras for some versions of the canonical anti-commutation
relations.

Recall that the CAR algebra with $d$ degrees of freedom is generated %
by $a_i,\ a_i^*$, $ i=1,\dots d$, and the relations
\begin{gather}
a_i^*a_i+a_ia_i^*=1,\ a_i^2=0,\ i=1,\dots d,\nonumber\\
a_i^*a_j=-a_ja_i^*,\ a_ja_i=-a_ia_j,\ i\ne j. \label{2}
\end{gather}
It is known that the Fock representation is the unique irreducible
representation of (\ref{2})
and the $C^*$-algebra generated by (\ref{2}) is isomorphic to $M_{2^d}(\mathbb{C})$.%

We  consider an interpolation between CCR and CAR known as $q$-CCR,
proposed by A.J.Macfarlane and L.C.Biedenharn for $d=1$, see
\cite{bied,mac} and by O.Greenberg, D.Fivel, M.Bozeiko and
R.Speicher for general $d$, see \cite{bsp,fiv,green}. Namely, the
higher-dimensional $q$-CCR have the following form
\begin{gather}\label{025}
a_i^*a_j=\delta_{ij}1+qa_ja_i^*, i=1,\dots,d,\ q\in(-1,1).
\end{gather}

%
Another well-known deformation of CAR, called twisted CAR, was
introduced and studied by W.Pusz, see \cite{p}. The twisted CAR
$*$-algebra (TCAR) is generated by $a_i,a_i^*,\ i=1,\dots,d$,
subject to the following relations
\begin{gather}
a_i^*a_i=1-a_i a_i^*-(1-\mu^2)\sum_{j<i}a_j a_j^*,\ i=1,\ldots,d,\ \mu\in(0,1),\nonumber\\
a_i^*a_j=-\mu a_j a_i^*,\ a_j a_i =-\mu a_i a_j,\ i<j,\label{055}\\
a_i^2=0,\ i=1,\ldots,d.\nonumber
\end{gather}

The Fock representation is the unique irreducible representation of
TCAR, and as in the non-deformed case the $C^*$-algebra generated by
TCAR coincides with $M_{2^d}(\mathbb{C})$.

The Wick analogue of TCAR (denoted below as WTCAR) one obtains from
TCAR taking away the relations between $a_i$, $a_j$. This algebra
was studied in \cite{jsw, pros}. In particular, it was shown that in
any representation of WTCAR the relations
\[
a_ja_i=-\mu a_ia_j,\ i<j,\ a_i^2=0,\ i=1,\dots,d-1,
\]
are satisfied and the irreducible representations of WTCAR were
classified.

The $q$-CCR and the WTCAR with $d=1$ are closely related with a
$*$-algebra known as the quantum disk. That is a $*$-algebra
generated by $a$ and $a^*$ satisfying the relation
\begin{gather}\label{09}
a^*a-qaa^ *=(1-q),\ q\in(-1,1).
\end{gather}
The family of $C^*$-algebras $D_q$ generated by ($\ref{09}$) was
studied by many authors, see for example, \cite{nn}.

If in (\ref{025}) we put $q=-1$ and in (\ref{055}) put $\mu=1$ we
get the Wick analogue of CAR, i.e. the $*$-algebra generated by
relations of the form
\begin{gather}\label{3}
a_i^*a_i+a_ia_i^*=1,\ i=1,\dots d,\\
a_i^*a_j=-a_ja_i^*,\ i\neq j.\nonumber
\end{gather}

In \cite{jw} P.E.T. J{\o}rgensen and R.F. Werner studied
representations of WCAR using the representation theory of Clifford
algebras. In particular, it was shown that in the irreducible
representations of (\ref{3}), for any pair $(i,j)$, one has
\[
a_ia_j+a_ja_i=y_{ij},\ \mbox{where}\ y_{ij}\in\mathbb{C},\
\norm{(y_{ij})}\leq 1.
\]
It was also stated in \cite {jw} that for
any $Y=(y_{ij}),\ y_{ji}=y_{ij}\in\mathbb{C}$ with $\norm{Y}\leq 1$
there exists an irreducible representation of (\ref{3}) with
\begin{gather}\label{4}
a_ia_j+a_ja_i=y_{ij}.
\end{gather}\noindent
Moreover,  the $C^*$-algebra $\widehat{\mathcal{E}}(-1,Y)$%
\quad generated by relations (\ref{3}), (\ref{4}) was shown to be
either isomorphic to $M_n(\mathbb{C})$ or to $M_n(\mathbb{C})\oplus
M_n(\mathbb{C})$ for appropriate $n\in \mathbb{N}$. In particular,
this fact implies that for any fixed $Y,\ \norm{Y}\leq 1,\ Y^T=Y$,
there are at most two non-equivalent irreducible representations of
relations (\ref{3}), (\ref{4}). Here we continue the study of
$C^*$-algebras associated with WCAR.

Our paper is organized as follows. In Section 1 we give some
definitions and facts used in the paper and fix notation.

In Section 2 we obtain a realization of WTCAR algebra as algebra of
continuous matrix-functions satisfying some boundary conditions.
Analysis of the case $d=1$ is crucial. Note that for $d=1$ we get,
up to normalization, the quantum disk with $q=-1$, called also the
``non-commutative circle''. The $C^*$-algebra $D_{-1}$ was studied
in \cite{nn}. In particular it was shown that $D_{-1}$ can be
faithfully  embedded into the algebra, $C(D^2\to M_2(\mathbb{C}))$,
of continuous matrix-functions on the unit disk. We make this result
more precise and show that $D_{-1}$ is isomorphic to an algebra of
continuous $2\times 2$ matrix-functions on the disk $D^2$ satisfying
certain boundary conditions on $S^1=\partial D^2$. Note that it is
more convenient for us to use the embedding of $D_{-1}$ into
$M_2(C(D^2))$ in the form different from the one presented in
\cite{nn}.

In the Section 3 we study representations of WTCAR with $d=2$ using
a dynamical systems technique, see \cite{osam}. For any
$y\in\mathbb{C},\ |y|\leq 1$ we give a parameterization of the
unitary equivalence classes of irreducible representations and
describe the  $C^*$-algebra $\mathcal{E}(-1,y)$ generated by
relations
\begin{gather}
a_i^*a_i+a_ia_i^*=1,\ i=1,2,\nonumber\\
a_1^*a_2=-a_2a_1^*, \label{wcr}\\
a_2a_1+a_1a_2=y\nonumber.
\end{gather}
Further we prove that the set of the isomorphism classes of
$\mathcal{E}(-1,y)$ consists of three elements:
$[\mathcal{E}(-1,0)]$, $[\mathcal{E}(-1,y),\ 0<|y|<1]$ and
$[\mathcal{E}(-1,y),\ |y|=1]$, where by $[\cdot]$ we denote the
class of isomorphic algebras. We also describe the $C^*$-algebras
$\mathcal{E}(\varepsilon)$, $0<\varepsilon<1$, defined by
(\ref{wcr}) where $y$ takes any value  from the set
$\{\varepsilon\le |y|\le1\}$. The isomorphism question is also
discussed.

In Section 4 we describe the enveloping $C^*$-algebra of WCAR.

\section{Preliminaries}
\noindent In this section, for convenience of the reader, we fix
some notation and recall necessary definitions and facts used in the
paper.

Let $\mathcal{A}$ be a $*$-algebra, having at least one
representation. Then a pair $(A,\rho)$ of a $C^*$-algebra $A$ and a
homomorphism $\rho:\mathcal{A}\rightarrow A$ is called an {\it
enveloping pair} for $\mathcal{A}$ if every irreducible
representation $\pi:\mathcal{A}\rightarrow B(H)$ factors uniquely
through $A$, i.e. there is a unique irreducible representation
$\pi_1$ of the algebra $A$ satisfying $\pi_1\circ\rho=\pi$. The
$C^*$-algebra $A$ is called an {\it enveloping} for $\mathcal{A}$.
An enveloping $C^*$-algebra for a $*$-algebra $\mathcal{A}$ is
unique and exists iff the set of bounded representations of
$\mathcal{A}$ is not empty and $\mathcal{A}$ is $*$-bounded, i.e.
for any $a\in \mathcal{A}$ one can find $C_a>0$ such that for any
bounded representation, $\pi$, $\Vert \pi(a)\Vert\le C_a$.

The following statement is a simple corollary of the non-commutative
analogue of the Stone-Weierstrass theorem (see \cite{f,v}).
\begin{theorem}\label{gsw}
Let $Y$ be a compact Hausdorff space. Let $C\subseteq B$ be
subalgebras of $A=C(Y\rightarrow M_n(\mathbb{C}))$. For every pair
$x_1, x_2\in Y$ define
$B(x_1,x_2)$ $(C(x_1,x_2)\ \mbox{respectively})$ as :%
\[
B(x_1,x_2):=\{(f(x_1),f(x_2))\in M_n(\mathbb{C})\times
M_n(\mathbb{C})\ |f\ \in B\},
\]
Then
\[
B=C\Leftrightarrow B(y_1,y_2)=C(y_1,y_2),\forall y_1,y_2\in Y.
\]
\end{theorem}

For representations $\pi_1$, $\pi_2$ of $*$-algebra $\mathcal{A}$ on
Hilbert spaces $\mathcal{H}(\pi_1)$ and $\mathcal{H}(\pi_2)$
respectively, let $C(\pi_1,\pi_2)$ be the space of intertwining
operators
\[
C(\pi_1, \pi_2)=\set{c\in
B(\mathcal{H}(\pi_2),\mathcal{H}(\pi_1)):\pi_1(a) c=c\pi_2(a),\
a\in\mathcal{A}}.
\]
Note that $C(\pi_1,\pi_2)=\set{0}$ iff $\pi_1$, $\pi_2$ are
disjoint, i.e. $\pi_1$, $\pi_2$ do not have unitary equivalent
subrepresentations.

For a $*$-algebra $\mathcal{A}\subset B(\mathcal{H})$ we denote by
$\mathcal{A}^\prime$ its commutant, i.e.
\[
\mathcal{A}^\prime=\{c\in B(\mathcal{H}):\ ca=ac,\
a\in\mathcal{A}\}.
\]
In what follows we will identify $\mathbb{C}^n\otimes\mathbb{C}^m$
with $\mathbb{C}^{nm}$ in such way that for $A=(a_{ij})\in
M_n(\mathbb{C})$ and $B\in M_m(\mathbb{C})$ the matrix $A\otimes B$
is equal to $(a_{ij}B)\in M_{nm}(\mathbb{C})$.

\section{Enveloping $C^*$-algebra for WTCAR}
In this Section we give a realization of enveloping $C^*$-algebra
for the $*$-algebra, $\mathcal{B}_{\mu}^{(d)}$, generated by WTCAR
as algebra of continuous matrix-functions.
\subsection{} We first study the case $d=1$. Evidently
\begin{gather}\label{rel}
\mathcal{B}_{\mu}^{(1)}=\mathbb{C}\langle a,a^*\ |\ a^*a+aa^*=1\rangle.
\end{gather}

Above we noted that this $C^*$-algebra is isomorphic to the
non-commutative circle $D_{-1}$ studied, in particular in \cite{nn}.
To get a realization of $D_{-1}$ as continuous matrix-functions we
use a classification of its irreducible representations up to
unitary equivalence. We use method of dynamical systems presented in
\cite{osam} in order to obtain this classification. Let $\pi$ be a
representation of (\ref{rel}). We consider the polar decomposition
of $\pi(a)=uc$. If $\pi$ is irreducible then (\ref{rel}) implies
that $\sigma(c^2)=\{x,1-x\}$, $0\le x\le\frac{1}{2}$, and
$u^2=e^{i\phi}\mathbf{1}$ if $x\ne 0$ and $u^2=0$, if \mbox{$x=0$}.
Moreover the eigenvalues of $c^2$ should have the same
multiplicities (see \cite{osam}) implying that the irreducible
representations with $\sigma(c^2)\ne\{\frac{1}{2}\}$ are
two-dimensional and the  irreducible representations with
$\sigma(c^2)=\{\frac{1}{2}\}$ are one-dimensional. Finally we have
the following list of irreducible representations:
\begin{itemize}
\item 2-dimensional:
\begin{gather}\label{24}
\pi_{x,\phi}(a)=\left(%
\begin{array}{cc}
  0 & e^{i\phi}\sqrt{x} \\
  \sqrt{1-x} & 0 \\
\end{array}%
\right),x\in [0,1/2),\phi\in[0,2\pi),
\end{gather}
\item 1-dimensional:
\begin{gather}\label{25}
\rho_\phi(a)=\frac{e^{i\phi}}{\sqrt{ 2}},\ \phi\in[0,2\pi).
\end{gather}
\end{itemize}
An alternative description can be found in \cite{mar,nn}.
\begin{remark}
Let
\begin{gather}\label{V}
V(\phi)=\frac 1 {\sqrt{2}} \left(%
\begin{array}{cc}
  e^{i\frac \phi 2} & -e^{i\frac \phi 2} \\
  1 & 1 \\
\end{array}%
\right).
\end{gather}
Then
\begin{gather*}
V^*(\phi)\pi_{\frac 1 2,\phi}(a)V(\phi)=\frac 1{\sqrt 2}%
\left(%
\begin{array}{cc}
  e^{i\frac{\phi}2} & 0 \\
  0 & -e^{\frac{i\phi}2} \\
\end{array}%
\right).
\end{gather*}
\end{remark}\noindent
showing that  any $1$-dimensional representation can be obtained
decomposing the representation $\pi_{\frac{1}2,\phi}$ with some
fixed $\phi$ into irreducible ones.

A result similar to one given in the next theorem can be found in
\cite{mar}. Here we anyway give a detailed proof of the statement,
since it presents in the most transparent way an idea of the more
tedious proofs of Theorems
\ref{theorem2}, \ref{theorem3}, \ref{theorem4}.
\begin{theorem}\label{theorem1}
The $C^*$-algebra $D_{-1}$ is isomorphic to the $C^*$-algebra
\begin{gather*}
A_1=\{f\in C(D^2\to M_2(\mathbb{C}))\ |\
V^*(\phi)f(e^{i\phi})V(\phi)\in\mathbb{C}\oplus\mathbb{C},
\forall\phi\in[0,2\pi]\},
\end{gather*}
where $D^2=\set{z\in\mathbb{C}:\ |z|\leq 1}$, and $V(\phi)$ is given
by the (\ref{V}).
\end{theorem}
\begin{proof} Let $I_{\frac{1}2}=[0,1/2]$ and $S^1=\set{z\in \mathbb{C}\ |\
|z|=1}$. Firstly we give a "primary" realization of $D_{-1}$. We
show, that $D_{-1}$ is isomorphic to
\begin{gather*}
A_0=\{f\in C(I_{\frac{1}2}\times S^1\to
M_2(\mathbb{C}))\ |\ f(0,e^{i\phi})=f(0,1),\ \forall\phi\in[0,2\pi],\\
V^*(\phi)f(1/2,e^{i\phi})V(\phi)\mbox{ is diagonal}, \forall\phi\in
[0,2\pi]\},
\end{gather*}

Let $\widetilde{a}=\widetilde{a}(x,e^{i\phi}):I_{\frac{1}2}\times
S^1\to M_2(\mathbb{C})$ be the function given by (\ref{24}), i.e.
\begin{gather*}
\widetilde{a}(x,e^{i\phi})=\pi_{x,\phi}(a)=\left(%
\begin{array}{cc}
  0 & e^{i\phi}\sqrt{x} \\
  \sqrt{1-x} & 0 \\
\end{array}%
\right).
\end{gather*}
One can check that $\widetilde{a}\in A_0$. Let $\widehat{A_0}$ be
the $C^*$-subalgebra of $A_0$ generated by $\widetilde{a}$. The
isomorphism $\widehat{A_0}\simeq D_{-1}$ follows directly from the
definition of enveloping pair.

To prove the equality $\widehat{A_0}=A_0$ we check the conditions of
Theorem \ref{gsw}, i.e.
\begin{gather}\label{27}
\widehat{A}_0((x_1,e^{i\phi
_1}),(x_2,e^{i\phi_2}))=A_0((x_1,e^{i\phi_1}),(x_2,e^{i\phi_2})),\\
\forall (x_1,e^{i\phi_1}),(x_2,e^{i\phi_2})\in I_{\frac{1}2}\times S^1.\nonumber%
\end{gather}
Since $\widehat{A}_0((x_1,e^{i\phi_1}),(x_2,e^{i\phi_2}))\subset
A_0((x_1,e^{i\phi_1}),(x_2,e^{i\phi_2}))$ and these algebras are
finite-dimensional, to prove (\ref{27}) it is sufficient to show
that their commutants are equal.

On the set $I_{\frac{1}2}\times S^1$ we introduce the equivalence
\begin{equation}\label{equiv}
(x_1,e^{i\phi_1})\sim(x_2,e^{i\phi_2}),\  \mbox{iff}\ \  x_1=x_2=0.
\end{equation}
Note, that if $(x_1,e^{i\phi_1})\nsim(x_2,e^{i\phi_2})$, then
\[
C(\pi_{x_1,\phi_1},\pi_{x_2,\phi_2})=\set{0}.
\]
Since
\[
\widehat{A}_0((x_1,e^{i\phi_1}),(x_2,e^{i\phi_2}))=\set{\left(%
\begin{array}{cc}
  \pi_{x_1,\phi_1}(b) & 0 \\
  0 & \pi_{x_2,\phi_2}(b) \\
\end{array}%
\right),\ b\in D_{-1}},
\]
then  %
\[
{\widehat{A}_0}^{\prime}((x_1,e^{i\phi_1}),(x_2,e^{i\phi_2}))=\set{\left(%
\begin{array}{cc}
  \Lambda_1 & 0 \\
  0 & \Lambda_2 \\
\end{array}%
\right),\ \Lambda_i\in\set{\pi_{x_i,\phi_i}(b),\ b\in D_{-1}}'}.
\]
The inclusion $\widehat{A}_0\subset A_0$ implies that
\[
A_0'((x_1,e^{i\phi_1}),(x_2,e^{i\phi_2}))=\set{\left(%
\begin{array}{cc}
  \Lambda_1 & 0 \\
  0 & \Lambda_2 \\
\end{array}%
\right),\ \Lambda_i\in A'_0(x_i,e^{i\phi_i})},
\]
where
\[
A_0(x,e^{i\phi})=\{f(x,e^{i\phi})\mid\ f\in A_0\}.
\]
If $(x_1,e^{i\phi_1})\sim(x_2,e^{i\phi_2})$ then $\ x_1=x_2=0$ and
$\pi_{x_1,\phi_1}=\pi_{x_2,\phi_2}=\pi_{0,1}$. In this case
\[
\widehat{A}_0'((0,1),(0,1))=\set{(\Lambda_{ij}),\ \Lambda_{ij}\in
\set{\pi_{0,1}(b)\ |\ b\in D_{-1}}',\ i,j=1,2}
\]
and
\[
A'_0((0,1),(0,1))=\set{(\Lambda_{ij}),\ \Lambda_{ij}\in A'_0(0,1),\
i,j=1,2}.
\]
It is left to show that for any $(x,\phi)\in
I_{\frac{1}{2}}\times[0,2\pi)$
\[
\set{\pi_{x,\phi}(b),\ b\in D_{-1}}'=A'_0(x,e^{i\phi}).%
\]
We consider two cases: $x\ne\frac{1}{2}$ and $x=\frac{1}{2}$.\\
\noindent \textbf{1)} Let $x\ne\frac{1}{2}$. Then $\pi_{x,\phi}$ is
irreducible and $\set{\pi_{x,\phi}(b),\ b\in
D_{-1}}=M_2(\mathbb{C})$ implying that
\[
\set{\pi_{x,\phi}(b),\ b\in D_{-1}}'=%
\set{\left(%
\begin{array}{cc}
  \lambda & 0 \\
  0 & \lambda \\
\end{array}%
\right),\ \lambda\in\mathbb{C}}=A'_0(x,\phi).
\]
\textbf{2)} Let $x=\frac{1}2$.

\textbf{2a)} If $\phi=0$ then $\pi_{\frac 1 2,0}(a)=\frac1{\sqrt 2}\left(%
\begin{array}{cc}
  0 & 1 \\
  1 & 0 \\
\end{array}%
\right)$. The boundary conditions for $f\in A_0$ imply that
$f(\frac 1 2,0)=\left(%
\begin{array}{cc}
   a & b \\
  b & a \\
\end{array}%
\right)$ for some $a,b\in\mathbb{C}$. Therefore,
\[
\set{\pi_{\frac{1}2,0}(D_{-1})}'=A'_0(1/2,0).
\]

\textbf{2b)} Let $\phi\in(0,2\pi),\ V=V(\phi)$ defined by (\ref{V}).
Then
\begin{gather*}
C\in\set{\pi_{\frac 1 2,\phi}(D_{-1})}'\Leftrightarrow
V^*CV\in\set{V^*\pi_{\frac 1 2,\phi}(D_{-1})V}'\Leftrightarrow\\
\Leftrightarrow V^*CV\in\set{\left(%
\begin{array}{cc}
  \lambda & 0 \\
  0 & \mu \\
\end{array}%
\right),\ \lambda,\mu\in\mathbb{C}}.
\end{gather*}
Analogously,
\[
C\in A'_0((1/2,\phi))\Leftrightarrow V^*CV\in\set{\left(%
\begin{array}{cc}
  \lambda & 0 \\
  0 & \mu \\
\end{array}%
\right),\ \lambda,\mu\in\mathbb{C}}.
\]
So we have $\widehat{A_0}=A_0$, and therefore $D_{-1}$ is isomorphic
to $A_0$. The quotient map by the equivalence (\ref{equiv}) induces
the isomorphism $A_0\simeq A_1$.

The proof is completed.
\end{proof}
\begin{remark}
From Theorem \ref{theorem1} one can get also a description of
$P(D_{-1})$, the dual space of $D_{-1}$ (compare with \cite[Theorem
1.1.]{nn})

Indeed, the isomorphism $D_{-1}\simeq A_0$ shows that the
$P_2(D_{-1})$, i.e the space of all pairwise non-equivalent
2-dimensional irreducible representations of $D_{-1}$, is
homeomorphic to the open disk $D^2\backslash\partial D^2$.

The injective map from the space of pairs
$(\rho_{\frac \phi 2},\rho_{\frac \phi 2+\pi})$ of 1-dimensional representations
into the space of all 2 dimensional representations, $T_2(D_{-1})$,
induces a covering over the circle $S^1=\partial D^2$. Our topology
on the dual space $P(D_{-1})$ can be described as follows. Let
$\mathcal{I}_2$ denote this covering (over $S^1$ with the structure
group and fiber $\mathbb{Z}_2$). Then $P(D_{-1})=\mathcal{I}_2\sqcup
(D^2\backslash\partial D^2)$, and a neighborhood of every $x\in
\mathcal{I}_2$ is the same as for $p(x)\in S^1=\partial D^2$ (it
implies that only points from the same fiber are non-separable, i.e.
for any such point there is no neighborhood which does not contain
the other point).

So, to define the topology on the $P(D_{-1})$ we have only to
determine the class of isomorphism of the covering $\mathcal{I}_2$.
From the formulas for representations one can see, that the total
space of covering is homeomorphic to $S^1$, so it is the unique
non-trivial one (the trivial coincides with $S^1\sqcup S^1$).
\end{remark}
\subsection{} Let us consider the case of general $d$. Using an
algorithm described in \cite{pros} one gets the following
\begin{proposition}\label{rpwtcar}
Any irreducible representation of WTCAR is unitarily equivalent to
one of the following
\\ \vspace{0.2 cm} \noindent
$\bullet\quad$ $2^d$-dimensional representations $\pi_{x,\phi}$,
$x\in[0,1/2),\ \phi\in[0,2\pi)$,
\begin{align}
\pi_{x,\phi}(a_i)=\bigotimes_{1\leq j<i}&\left(%
\begin{array}{cc}
   1 & 0 \\
  0 & -\mu \\
\end{array}%
\right)\otimes
\left(%
\begin{array}{cc}
   0 & 0 \\
  1 & 0 \\
\end{array}%
\right)\otimes %
\bigotimes_{i<j\leq d}
 \left(%
\begin{array}{cc}
   1 & 0 \\
  0 & 1 \\
\end{array}%
\right),\ i=1,\dots,d-1,\nonumber\\
\ & \label{repwtcar}\\
\pi_{x,\phi}(a_d)& =\bigotimes_{1\leq j<d}\left(%
\begin{array}{cc}
   1 & 0 \\
  0 & -\mu \\
\end{array}%
\right)\otimes
\left(%
\begin{array}{cc}
   0 & e^{i\phi}\sqrt{x} \\
  \sqrt{1-x} & 0 \\
\end{array}%
\right). \nonumber
\end{align}
\vspace{0.1cm}\noindent $\bullet\quad$ $2^{d-1}$-dimensional
representations $\rho_{\phi}$, $\phi\in[0,2\pi)$,
\begin{align}
\rho_\phi(a_i)=\bigotimes_{1\leq j<i}\left(%
\begin{array}{cc}
  1 & 0 \\
  0 & -\mu \\
\end{array}%
\right)& \otimes
\left(%
\begin{array}{cc}
  0 & 0 \\
  1 & 0 \\
\end{array}%
\right)\otimes %
\bigotimes_{i<j<d}
\left(%
\begin{array}{cc}
  1 & 0 \\
  0 & 1 \\
\end{array}%
\right),\ i=1,\dots,d-1,\nonumber \\
\ & \label{repwtcar1}\\
\rho_\phi(a_d)&=\frac{e^{i\phi}}{\sqrt 2}
\bigotimes_{1\leq j<d}\left(%
\begin{array}{cc}
  1 & 0 \\
  0 & -\mu \\
\end{array}%
\right). \nonumber
\end{align}

Moreover for different pairs $(x,\phi),\ x\in(0,1/2)$, the
corresponding representations $\pi_{x,\phi}$ are non-equivalent;
$\pi_{0,\phi_1}\stackrel{u}{\sim}\pi_{0,\phi_2}$ for any
$\phi_1,\phi_2\in [0,2\pi)$; for different $\phi$ the
representations $\rho_{\phi}$ are non-equivalent.
\end{proposition}
\begin{remark}
Letting $x=\frac{1}{2}$ in (\ref{repwtcar}) we get also a
representation, $\pi_{\frac{1}{2},\phi}$, of \mbox{WTCAR} but a
reducible one:
\[
\pi_{\frac{1}{2},\phi}\stackrel{u}{\sim}\rho_{\frac \phi
2}\oplus\rho_{\frac \phi 2 +\pi}.
\]
\end{remark}
\begin{theorem}
The enveloping $C^*$-algebra $C^*(B_{\mu}^{(d)})$ for \mbox{WTCAR}
is isomorphic to \mbox{$ M_{2^{d-1}}(\mathbb{C})\otimes D_{-1}$}.
\end{theorem}
\begin{proof}
As for the case $d=1$, the enveloping $C^*$-algebra coincides with
the $C^*$-algebra generated by matrix-functions on
$I_{\frac{1}{2}}\times S^1$  given by (\ref{repwtcar}). From
(\ref{repwtcar}) it also follows that
\begin{gather*}
C^*(\pi_{x,\phi}(a_i),\ i=1\ldots d-1)=\bigotimes_{i=1}^{d-1}
M_2(\mathbb{C})\otimes \mathbf{1},
\end{gather*}
and
\[
\bigotimes_{i=1}^{d-1}\mathbf{1}\otimes\left(%
\begin{array}{cc}
   0 & e^{i\phi}\sqrt{x} \\
  \sqrt{1-x} & 0 \\
\end{array}%
\right)\in C^*(\pi_{x,\phi}(a_i),\ i=1\ldots
d)=C^*(\mathcal{B}_{\mu}^{(d)})
\]
giving $C^*(\mathcal{B}_{\mu}^{(d)})=M_{2^{d-1}}(\mathbb{C})\otimes
D_{-1}$.
\end{proof}
It follows from Theorem 2 that the dual space for
$C^*(\mathcal{B}_{\mu}^{(d)})$ is the same as for the algebra
$D_{-1}$. We have also
\begin{corollary}
$C^*(\mathcal{B}_{\mu}^{(d)})\simeq C^*(\mathcal{B}_{0}^{(d)})$,
$0<\mu<1$.
\end{corollary}
\section {Representations of WCAR with two degrees of freedom}
In this section we study representations of WCAR with $d=2$ and
describe the corresponding families of $C^*$-algebras. Let
$\mathcal{E}(-1)$ denote the enveloping $C^*$-algebra of WCAR for
$d=2$:
\begin{gather}\label{**}
\mathcal{E}(-1)=C^*
\langle %
a^*_ia_i+a_i a^*_i=1,i=1,2,\ a^*_1a_2=-a_2a^*_1%
\rangle.
\end{gather}
Note, that WCAR is $*$-bounded: $\Vert \pi(a_i)\Vert\le 1$ for each
representation $\pi$ and $i=1,2$.

Let us first describe irreducible representations of
$\mathcal{E}(-1)$. In what follows if $\pi$ is a representation of
WCAR, we write simply $a_i$ instead of  $\pi(a_i)$, when no
confusion can arise.

Let $A=a_2a_1+a_1a_2$, then $A^*A=AA^*$ and $a^*_iA=Aa^*_i,\ i=1,2$,
which allows us to apply Fuglede's theorem (see, e.g., \cite{rud})
and get $Aa_i=a_iA$, $i=1,2$. So in irreducible representation we
must have
\begin{gather}\label{ygryk}
 a_1a_2+a_2a_1=y
\end{gather}
for some $y\in\mathbb{C}$. Moreover we will show below that $|y|\le
1$.

In the sequel we denote by $\mathcal{E}(-1,y)$ the quotient
\[
\mathcal{E}(-1)/\langle a_1a_2+a_2a_1-y\rangle.
\]
\subsection{Representations of $\mathcal{E}(-1,y),\ y\neq 0$.}
It is easy to check that $a^2_1$ is normal and
$a^*_ia^2_1=a^2_1a^*_i,\ i=1,2$. Then by Fuglede's theorem,
$a_1^2a_i=a_ia_1^2$, and Schur's lemma implies that $a^2_1=\rho_1
e^{i\phi_1}$.

Note that in irreducible representation of $\mathcal{E}(-1)$ we
have either $a_1^2=(a_1^*)^2=0$ or $\ker a_1=\ker a_1^*=\set{0}$.
Indeed, $\ker a_1^2\neq \set{0}\Leftrightarrow \ker
a_1\neq\set{0}$, and $\ker a_1^2$ is invariant with respect to
$a_i,a_i^*,\ i=1,2$.

Let $a_1:=u_1c_1$ be the polar decomposition of the operator $a_1$.
Then re\-la\-tions (\ref{**}) imply that
\[
c^2_1u_1=u_1(1-c^2_1)\ \mbox{\text and}\ c_1u_1=u_1\sqrt{(1-c^2_1)}.
\]
Hence
\[
a^2_1=u_1c_1u_1c_1=u^2_1c_1\sqrt{(1-c^2_1)}=\rho e^{i\phi},
\]
with $\rho\neq 0$ if $a_1^2\neq 0$.\\
\ \\ \noindent \textbf{a)} If $a_1^2=0$ then $a_1^*a_1+a_1a_1^*=1$
implies that
\[
a_1=\left(%
\begin{array}{cc}
  0 & 0 \\
  1 & 0 \\
\end{array}%
\right)\otimes\mathbf{1}.
\]
\textbf{b)} If $a_1^2\neq 0$ then by uniqueness of the polar
decomposition
\[
u_1^2=e^{i\phi_1}\cdot\mathbf{1},\
c_1^2(1-c_1^2)=\rho_1^2\cdot\mathbf{1}
\]
implying that $\sigma(c_1^2)=\set{x_1,1-x_1}$, where
$x_1(1-x_1)=\rho_1^2$, $0<x_1\leq 1/2$.

Further, if $\sigma(c_1^2)=\set{\frac 1 2}$ then $a_1=\frac 1{\sqrt
2}u_1$, where $u_1$ is unitary, and by Fuglede's theorem
$a_2a_1+a_1a_2=0$. So when $y\neq 0$ we should consider only the
case \mbox{$0<x_1<\frac 1{2}$}. Then, up to unitary equivalence,
\begin{equation}\label{a1}
a_1=\left(%
\begin{array}{cc}
  0 & e^{i\phi_1}\sqrt{x_1}\cdot\mathbf{1} \\
  \sqrt{1-x_1}\cdot\mathbf{1} & 0 \\
\end{array}%
\right),\ \phi_1\in[0,2\pi).
\end{equation}
We start with the analysis of case \textbf{b)} and suppose that
$a_1$ is given by  (\ref{a1})

From $a_1^*a_2=-a_2a_1^*$ we get
\[
a_2=\left(%
\begin{array}{cc}
  b_1 & b_2\sqrt{1-x_1} \\
  -e^{-i\phi_1}b_2\sqrt{x_1} & -b_1 \\
\end{array}%
\right),
\]
then the relation $a_1a_2+a_2a_1=y$ implies that $b_2(1-2x_1)=y$.
Since $0<x_1<\frac 1 2$,
\[
b_2=\frac{y}{1-2x_1},\ a_2=\left(%
\begin{array}{cc}
  b_1 & \frac{y}{1-2x_1}\sqrt{1-x_1} \\
  -e^{-i\phi_1} \frac{y}{1-2x_1}\sqrt x_1 & -b_1 \\
\end{array}%
\right),
\]
and from $a_2^*a_2+a_2a_2^*=1$ we get
\begin{gather}\label{bb}
b_1^*b_1+b_1b_1^*=1-\frac{|y|^2}{(1-2x_1)^2}.
\end{gather}
Clearly the family $\set{a_i,a_i^*,\ i=1,2}$ is irreducible iff so
is  $\set{b_1,b_1^*}$.

Representations of (\ref{bb}) exist iff  $|y|\leq 1-2x_1$.
Furthermore $x_1\in(0,1/2)$ implies $|y|\leq 1$ and
$x_1\leq\frac{1-|y|}2$. Evidently, $0<x_1\leq\frac{1-|y|}2,\ |y|\leq
1$ yields $|y|<1$.\\ \noindent
\textbf{b1)} If $x_1=\frac{1-|y|}2$, then $b_1=0$.\\
\noindent \textbf{b2)} If $0<x_1<\frac{1-|y|}2$ put
$b_1=\sqrt{1-\frac{|y|^2}{(1-2x_1)^2}}\cdot\widetilde{b}_1$, then
$\widetilde{b}^*_1\widetilde{b}_1+\widetilde{b}_1\widetilde{b}^*_1=1$
and the family $\set{a_i,\ a_i^*,\ i=1,2}$ is irreducible iff so is
$\set{\widetilde{b}_1,\ \widetilde{b}^*_1}$. Two families
$\set{a_i^{(j)},\ {a_i^{(j)*}},\ i=1,2}$, $j=1,2$, are unitarily
equivalent iff the corresponding families
$\set{\widetilde{b}_1^{(j)},\ {{\widetilde{b}_1}^{(j)*}}}$, \mbox{$j=1,2$},
are unitarily equivalent.

Then in the irreducible case one has (see Section 2)

\textbf{b2.1)}
\begin{gather}b_1=\sqrt{1-\frac{|y|^2}{(1-2x_1)^2}}\left(%
\begin{array}{cc}
  0 & e^{i\phi_2}\sqrt {x_2} \\
  \sqrt{1-x_2} & 0 \\
\end{array}%
\right),\label{f1}\\
0\leq x_2<\frac{1}2,\ \phi_2\in[0,2\pi).\nonumber
\end{gather}

\textbf{b2.2)}
\begin{equation}\label{f2}
b_1=\frac 1{\sqrt 2}\sqrt{1-\frac{|y|^2}{(1-2x_1)^2}}e^{i\phi_2},\
\phi_2\in[0,2\pi).
\end{equation}
We now turn to case \textbf{a)}: $x_1=0$. Repeating the arguments
given in {\bf b)} with $x_1=0$ we obtain that
\[
a_1=\left(%
\begin{array}{cc}
  0 & 0 \\
  \mathbf{1} & 0 \\
\end{array}%
\right),\ %
a_2=\left(%
\begin{array}{cc}
  b_1 & y \\
  0 & -b_1 \\
\end{array}%
\right)\mbox{ and } b_1^*b_1+b_1b_1^*=1-|y|^2.
\]
If $|y|=1$, then $b_1=0$. When $|y|<1$ to get $b_1$ one has to put
$x_1=0$ in (\ref{f1}), (\ref{f2}).

Let
\[
D_y=\set{(x_1,x_2)|\ 0\leq x_1\leq\frac{1-|y|} 2,\ %
 0\leq x_2\leq\frac{1}2}
\]
Define an equivalence on $D_y \times\mathbf{T}^2$:
\begin{equation*}
(x_1,x_2,e^{i\phi_1},e^{i\phi_2})\sim(y_1,y_2,e^{\psi_1},e^{i\psi_2})
\end{equation*}
iff either $x_1=y_1=0$, $x_2=y_2$, $\phi_2=\psi_2$ or $x_2=y_2=0$,
$x_1=y_1$, $\phi_1=\psi_1$ or $x_1=y_1=\frac{1-|y|}{2}$,
$\phi_1=\psi_1$.

Summarizing the above discussion we have
\begin{theorem}
The $C^*$-algebra $\mathcal{E}(-1,y)$ is non-zero iff $|y|\le 1$. If
$y\ne 0$, then any irreducible representation of $\mathcal{E}(-1,y)$
is unitarily equivalent to one of the following:\\ \noindent
\textbf{a)} If $|y|=1$, $
a_1=\left(%
\begin{array}{cc}
  0 & 0 \\
  1 & 0 \\
\end{array}%
\right),\ %
a_2=\left(%
\begin{array}{cc}
  0 & y \\
  0 & 0 \\
\end{array}%
\right). $\\ \noindent \textbf{b)} $|y|<1$.
\vspace{0.2cm}\\
\textbf{b1)} If $x_1=\frac{1-|y|}2$, then
\begin{gather*}
a_1=\left(%
\begin{array}{cc}
  0 & e^{i\phi_1}\sqrt{\frac{1-|y|}2} \\
  \sqrt{\frac{1+|y|}2} & 0 \\
\end{array}%
\right),\ %
a_2=\left(%
\begin{array}{cc}
  0 & \frac y{|y|}\sqrt{\frac{1+|y|}2} \\
  -\frac {e^{-i\phi_1}y}{|y|}\sqrt{\frac{1-|y|}2} & 0 \\
\end{array}%
\right).
\end{gather*}
\textbf{b2)} If $0\leq x_1<\frac{1-|y|}2$, then
\vspace{0.2cm}\\

\textbf{b2.1)} if $0\leq x_2<\frac{1}2$,
\begin{align}\label{cch}
a_1 & =\left(%
\begin{array}{cc}
  0 & e^{i\phi_1}\sqrt{x_1} \\
  \sqrt{1-x_1} & 0\\
\end{array}%
\right)\otimes%
\left(%
\begin{array}{cc}
  1 & 0 \\
  0 & 1 \\
\end{array}%
\right),\ \phi_1\in [0,2\pi),\ \phi_2\in [0,2\pi)\nonumber\\  %
\  \nonumber \\
a_2 & =\sqrt{1-\frac{|y|^2}{(1-2x_1)^2}}%
\left(%
\begin{array}{cc}
  1 & 0 \\
  0 & -1 \\
\end{array}%
\right)\otimes
\left(%
\begin{array}{cc}
  0 & e^{i\phi_2}\sqrt{x_2} \\
  \sqrt{1-x_2} & 0 \\
\end{array}%
\right)+  \\
\ \nonumber \\
&+\frac{y}{1-2x_1}\left(%
\begin{array}{cc}
  0 & \sqrt{1-x_1} \\
  -e^{-i\phi_1}\sqrt{x_1} & 0 \\
\end{array}%
\right)\otimes\left(%
\begin{array}{cc}
  1 & 0 \\
  0 & 1 \\
\end{array}%
\right);\nonumber
\end{align}

\textbf{b2.2)} if $x_2=\frac{1}2$,
\begin{gather*}
a_1=\left(%
\begin{array}{cc}
  0 & e^{i\phi_1}\sqrt{x_1} \\
  \sqrt{1-x_1} & 0\\
\end{array}%
\right),\\
\ \\
a_2=\left(%
\begin{array}{cc}
  \frac {e^{i\phi_2}}{\sqrt 2}\sqrt{1-\frac{|y|^2}{(1-2x_1)^2}} & \frac{y\sqrt{1-x_1}}{1-2x_1} \\
  -e^{-i\phi_1}\frac{y\sqrt{x_1}}{1-2x_1} & -\frac {e^{i\phi_2}}{\sqrt 2}\sqrt{1-\frac{|y|^2}{(1-2x_1)^2}} \\
\end{array}%
\right),\ \phi_1,\phi_2\in[0,2\pi).
\end{gather*}
When $y\ne 0,\ |y|<1$, is fixed, the representations corresponding
to non-equivalent tuples $(x_1,x_2,e^{i\phi_1},e^{i\phi_2})\in D_y
\times\mathbf{T}^2$ are non-equivalent.
\end{theorem}
Using this classification we can describe $\mathcal{E}(-1,y)$ for
 $0<|y|\le 1$. The case $y=0$ will be
studied separately below.
\begin{theorem}\label{theorem2}
\textbf{1)} If $0<|y|<1$ then $\mathcal{E}(-1,y)$ can be realized as
follows
\begin{gather*}
\mathcal{E}(-1,y)=\bigl\{f\in C(D_y\times\mathbf{T}^2\to M_4(\mathbb{C}))\ |\\
f(0,x_2,e^{i\phi_1},e^{i\phi_2})=f(0,x_2,1,e^{i\phi_2}),\\
f(x_1,0,e^{i\phi_1},e^{i\phi_2})=f(x_1,0,e^{i\phi_1},1),\\
f\bigl((1-\mid y\mid)/2,x_2,e^{i\phi_1},e^{i\phi_2}\bigr)=
f\bigl((1-\mid y\mid)/2,0,e^{i\phi_1},1\bigr)\in M_2(\mathbb{C})\otimes\mathbf{1},\\
\nu^*_2(\phi_2)f(x_1,1/2,e^{i\phi_1},e^{i\phi_2})\nu_2(\phi_2)\in
M_2(\mathbb{C})\oplus M_2(\mathbb{C}),\ \forall\phi_2\in[0,2\pi]
\bigr\},
\end{gather*}
where
\begin{gather*}
\nu_2:[0,2\pi]\to U(4),\ \phi\mapsto T(V(\phi)\otimes\mathbf{1}_2 ),\\
T:\mathbb{C}^2\otimes\mathbb{C}^2\to\mathbb{C}^2\otimes\mathbb{C}^2,\ x\otimes y\mapsto y\otimes x, %
\ \mbox{ $V(\phi)$ is given by (\ref{V})}.
\end{gather*}
\textbf{2)} If $|y|=1$ then $\mathcal{E}(-1,y)$ is isomorphic to
$M_2(\mathbb{C})$.
\end{theorem}
\begin{proof}
First, note that any irreducible representation of
$\mathcal{E}(-1,y)$ is either defined by formulas (\ref{cch}) or can
be obtained by decomposition in irreducible components of such
representations with $x_1=\frac{1-\mid y\mid}{2}$ or
$x_2=\frac{1}{2}$. Hence if $a_i$ are given by (\ref{cch}),
$C^*\bigl(a_i(x_1,x_2,e^{i\phi_1},e^{i\phi_2}),\ i=1,2\bigr)$ is
isomorphic to $\mathcal{E}(-1,y)$. Further, it is easy to check that
\begin{gather*}
\nu_2^*(\phi_2)a_1(x_1,1/2,e^{i\phi_1},e^{i\phi_2})\nu_2(\phi_2)=\\
\ \\
=\left(%
\begin{array}{cc}
  1 & 0 \\
  0 & 1 \\
\end{array}%
\right)\otimes%
\left(%
\begin{array}{cc}
  0 & e^{i\phi_1}\sqrt{x_1} \\
  \sqrt{1-x_1} & 0 \\
\end{array}%
\right)\in M_2(\mathbb{C})\oplus M_2(\mathbb{C}),\\
\ \\
\nu_2^*(\phi_2)a_2(x_1,1/2,e^{i\phi_1},e^{i\phi_2})\nu_2(\phi_2)=\\%
\ \\
=\frac 1{\sqrt 2}%
\sqrt{1-\frac{|y|^2}{(1-2x_1)^2}}
\left(%
\begin{array}{cc}
  e^{i\frac{\phi_2}2} & 0 \\
  0 & -e^{i\frac{\phi_2}2} \\
\end{array}%
\right)\otimes
\left(%
\begin{array}{cc}
  1 & 0 \\
  0 & -1 \\
\end{array}%
\right)+\\
\ \\
+\frac y{1-2x_1}%
\left(%
\begin{array}{cc}
  1 & 0 \\
  0 & 1 \\
\end{array}%
\right)\otimes%
\left(%
\begin{array}{cc}
  0 & \sqrt{1-x_1} \\
  e^{-i\phi_1} & 0 \\
\end{array}%
\right)\in M_2(\mathbb{C})\oplus M_2(\mathbb{C}),\\
\ \\
a_1\bigl((1-\mid y\mid)/2,x_2,e^{i\phi_1},e^{i\phi_2}\bigr)=\left(%
\begin{array}{cc}
  0 & e^{i\phi_1}\sqrt{\frac{1-|y|}2} \\
  \sqrt{\frac{1+|y|}2} & 0 \\
\end{array}%
\right)\otimes\mathbf{1}_2,\\
\ \\
a_2\bigl((1-\mid y\mid)/2,x_2,e^{i\phi_1},e^{i\phi_2}\bigr)=\left(%
\begin{array}{cc}
  0 & \frac y{|y|}\sqrt{\frac{1+|y|}2} \\
  -\frac {e^{-i\phi_1}y}{|y|}\sqrt{\frac{1-|y|}2} & 0 \\
\end{array}%
\right)\otimes\mathbf{1}_2.
\end{gather*}
and if $\pi_j,\ j=1,2$, are representations (possibly reducible)
defined by (\ref{cch}) with non-equivalent tuples of parameters then
$C(\pi_1,\pi_2)=\set{0}$. To complete the proof use the same
arguments as in the proof of Theorem 1.
\end{proof}
\subsection{Representations and enveloping $C^*$-algebra of $\mathcal{E}(-1,0)$}
Let us consider the case $y=0$. Analysis similar to that in
Subsection 3.1 gives  a description of the unitary equivalence
classes of irreducible representations of $\mathcal{E}(-1,0)$. Note
that  here one must also consider the case
$\sigma(c_1^2)=\set{\frac{1}2}$. Namely, we have that the generators
$a_1,a_2$ in irreducible representations of $\mathcal{E}(-1,0)$ have
one of the forms presented below.\\ \noindent \textbf{1)}
\begin{gather}
a_1=\left(%
\begin{array}{cc}
  0 & e^{i\phi_1} \sqrt{x_1} \\
  \sqrt{1-x_1} & 0 \\
\end{array}%
\right)\otimes%
\left(%
\begin{array}{cc}
  1 & 0 \\
  0 & 1 \\
\end{array}%
\right),\nonumber\\
\ \label{regg1} \\ %
a_2=\left(%
\begin{array}{cc}
  1 & 0 \\
  0 & -1 \\
\end{array}%
\right)\nonumber
\otimes%
\left(%
\begin{array}{cc}
  0 & e^{i\phi_2} \sqrt{x_2} \\
  \sqrt{1-x_2} & 0 \\
\end{array}%
\right),\\
\ \nonumber\\
0\leq x_1,x_2 < 1/2,\ 0 \leq \phi_1,\phi_2 < 2\pi;\nonumber
\end{gather}
\textbf{2)}

\textbf{2a)}
\begin{gather}
a_1=\frac{e^{i\phi_1}}{\sqrt{2}}\left(%
\begin{array}{cc}
  1 & 0 \\
  0 & -1 \\
\end{array}%
\right),%
a_2=\left(%
\begin{array}{cc}
0 & e^{i\phi_2} \sqrt{x_2} \\
  \sqrt{1-x_2} & 0 \\
\end{array}%
\right),\label{2a}\\ \ \nonumber\\ 0\leq x_2 < 1/2,\ 0 \leq
\phi_1,\phi_2 < 2\pi,\nonumber
\end{gather}

\textbf{2b)}
\begin{gather*}
a_1=\left(%
\begin{array}{cc}
0 & e^{i\phi_1} \sqrt{x_1} \\
  \sqrt{1-x_1} & 0 \\
\end{array}%
\right),
a_2=\frac{e^{i\phi_2}}{\sqrt{2}}\left(%
\begin{array}{cc}
  1 & 0 \\
  0 & -1 \\
\end{array}%
\right),%
\\
\ \\
0\leq x_1 < 1/2,\ 0 \leq \phi_1,\phi_2 < 2\pi;
\end{gather*}
\textbf{3)}
\begin{gather}
a_1=\frac{e^{i\phi_1}}{\sqrt{2}}\left(%
\begin{array}{cc}
  1 & 0 \\
  0 & -1 \\
\end{array}%
\right),%
a_2=\frac1{\sqrt 2}\left(%
\begin{array}{cc}
  0 & e^{i\phi_2} \\
  1 & 0 \\
\end{array}%
\right),\label{31}\\ \ \nonumber \\ 0 \leq \phi_1 < \pi,\ 0\leq
\phi_2<2\pi.\nonumber
\end{gather}
As above we introduce an equivalence on $I_{\frac 1 2}\times
I_{\frac 1 2}\times\mathbf{T}^2$:
\[
(x_1,x_2,e^{i\phi_1},e^{i\phi_2})\sim
(y_1,y_2,e^{i\psi_2},e^{i\psi_2}),\
\]
iff $x_1=y_1=0$, $x_2=y_2$, $\phi_2=\psi_2$ or $x_2=y_2=0$,
$x_1=y_1$, $\phi_1=\psi_1$. The representations given by {\bf 1),
2), 3)} are equivalent iff they correspond to equivalent quadruples
of parameters.

We next present a realization of $\mathcal{E}(-1,0)$ analogous to
that given in \mbox{Section 3.1.}
\begin{theorem}\label{theorem3}
The $C^*$-algebra $\mathcal{E}(-1,0)$ is isomorphic to the following
algebra of continuous matrix-functions
\begin{gather*}
\bigl\{f\in C(I_{\frac 1 2}\times I_{\frac 1 2}\times S^1\times S^1\to
M_4(\mathbb{C}))\ |\\
f(0,x_2,e^{i\phi_1},e^{i\phi_2})=f(0,x_2,1,e^{i\phi_2}),\\
f(x_1,0,e^{i\phi_1},e^{i\phi_2})=f(x_1,0,e^{i\phi_1},1),\\
\nu^*_1(\phi_1,\phi_2)f(1/2,x_2,e^{i\phi_1},
e^{i\phi_2})\nu_1(\phi_1,\phi_2)\in M_2(\mathbb{C})\oplus
M_2(\mathbb{C}),\ \forall \phi_1,\phi_2\in[0,2\pi],
\\
\nu^*_2(\phi_2)f(x_1,1/2,e^{i\phi_1},e^{i\phi_2})\nu_2(\phi_2)\in
M_2(\mathbb{C})\oplus M_2(\mathbb{C}), \forall
\phi_2\in[0,2\pi]\bigr\},
\end{gather*}
where
\begin{gather*}
\nu_1:[0,2\pi]\times[0,2\pi]\to U(4),\ (\phi_1,\phi_2)\mapsto W(\phi_2)(V(\phi_1) \otimes \mathbf{1}_2),\\
\nu_2:[0,2\pi)\to U(4),\ \phi_2\mapsto T(V(\phi_2)\otimes\mathbf{1}_2 ),\\
T:\mathbb{C}^2\otimes\mathbb{C}^2\to\mathbb{C}^2\otimes\mathbb{C}^2,\ x\otimes y\mapsto y\otimes x,\\
W(\phi)=(\mathbf{1}_2\otimes V(\phi))S(\mathbf{1}_2\otimes V^*(\phi)),\ \mbox{$V(\phi)$ is given by (\ref{V})},\\ %
\ \\
S=\left(%
\begin{array}{cc}
  1 & 0 \\
  0 & 0 \\
\end{array}%
\right)\otimes%
\left(%
\begin{array}{cc}
  1 & 0 \\
  0 & 1 \\
\end{array}%
\right)+%
\left(%
\begin{array}{cc}
  0 & 0 \\
  0 & 1 \\
\end{array}%
\right)\otimes%
\left(%
\begin{array}{cc}
  0 & 1 \\
  1 & 0 \\
\end{array}%
\right).
\end{gather*}
\end{theorem}\nopagebreak
\begin{proof}
The proof is similar to ones of Theorems \ref{theorem1} and
\ref{theorem2}. We note only that
\begin{gather*}
\nu^*_1(\phi_1,\phi_2)a_1(1/2,x_2,e^{i\phi_1},e^{i\phi_2})\nu_1(\phi_1,\phi_2)=\\
\ \\
=\frac 1{\sqrt 2}\left(%
\begin{array}{cc}
  e^{i\frac{\phi_1}2} & 0 \\
  0 & e^{-i\frac{\phi_1}2} \\
\end{array}%
\right)\otimes%
\left(%
\begin{array}{cc}
  1 & 0 \\
  0 & -1 \\
\end{array}%
\right)%
\in M_2(\mathbb{C})\oplus M_2(\mathbb{C}),\\
\ \\
\nu^*_1(\phi_1,\phi_2)a_2(1/2,x_2,e^{i\phi_1},e^{i\phi_2})\nu_1(\phi_1,\phi_2)=\\
\ \\
=\left(%
\begin{array}{cc}
  1 & 0 \\
  0 & 1 \\
\end{array}%
\right)\otimes%
\left(%
\begin{array}{cc}
  0 & e^{i\phi_2}\sqrt{x_2} \\
  \sqrt{1-x_2} & 0 \\
\end{array}%
\right) \in M_2(\mathbb{C})\oplus M_2(\mathbb{C}),
\end{gather*}
\begin{gather*}
\nu^*_2(\phi_2)a_1(x_1,1/2,e^{i\phi_1},e^{i\phi_2})\nu_2(\phi_2)=\\
\ \\
=\left(%
\begin{array}{cc}
  1 & 0 \\
  0 & 1 \\
\end{array}%
\right)\otimes%
\left(%
\begin{array}{cc}
  0 & e^{i\phi_1}\sqrt{x_1} \\
  \sqrt{1-x_1} & 0 \\
\end{array}%
\right)
\in M_2(\mathbb{C})\oplus M_2(\mathbb{C}),\\
\ \\
\nu^*_2(\phi_2)a_2(x_1,1/2,e^{i\phi_1},e^{i\phi_2})\nu_2(\phi_2)=\\
\ \\
=\frac 1{\sqrt 2}\left(%
\begin{array}{cc}
  e^{i\frac{\phi_2}2} & 0 \\
  0 & e^{-i\frac{\phi_2}2} \\
\end{array}%
\right)\otimes%
\left(%
\begin{array}{cc}
  1 & 0 \\
  0 & -1 \\
\end{array}%
\right)%
\in M_2(\mathbb{C})\oplus M_2(\mathbb{C}).
\end{gather*}
and $C(\pi_1,\pi_2)=\{0\}$ if $\pi_1$, $\pi_2$ are determined by
(\ref{regg1}) with non-equivalent tuples of parameters.
\end{proof}
Next we study when the $C^*$-algebras $\mathcal{E}(-1,y)$ are
isomorphic.
\begin{proposition}
For any $y_1,\ y_2$, such that either $0<|y_i|<1$, $i=1,2$, or
$|y_i|=1$, $i=1,2$, $\mathcal{E}(-1,y_1)$ is isomorphic to
$\mathcal{E}(-1,y_2)$.  $\mathcal{E}(-1,0)$ is not isomorphic to any
other $\mathcal{E}(-1,y)$.
\end{proposition}
\begin{proof}
The isomorphism of  $\mathcal{E}(-1,y_1)$ and $\mathcal{E}(-1,y_2)$
with  $0<|y_i|<1$, $i=1,2$, is induced by the homeomorphism
$\vartheta\colon D_{y_1}\times\mathbf{T}^2\rightarrow
D_{y_2}\times\mathbf{T}^2$ given by the rule
\[
(x_1,x_2,e^{i\phi_1},e^{i\phi_2})\mapsto
(\frac{1-|y_2|}{1-|y_1|}x_1,x_2,e^{i\phi_1},e^{i\phi_2})
\]
To see that $\mathcal{E}(-1,0)$ is not isomorphic to
$\mathcal{E}(-1,y)$ when $0<|y|<1$ consider the family
$\mathbf{\rho}=\{\rho_{\lambda_1,\lambda_2},\ |\lambda_i|=1,\
i=1,2\}$ of automorphisms of both algebras defined by
\[
\rho_{\lambda_1,\lambda_2}(a_i)=\lambda_i a_i,\ i=1,2.
\]
Denote by $\mathcal{F}_0\subset\mathcal{E}(-1,0)$ and
$\mathcal{F}_y\subset\mathcal{E}(-1,y)$ the $C^*$-subalgebras of the
elements fixed by the family $\mathbf{\rho}$.

Then
\[
\mathcal{F}_0=C^*(a_i^n (a_i^*)^n,\ (a_i^*)^n a_i^n,\ i=1,2,\
n\in\mathbb{Z}_{+}).
\]
In fact, given $b\in\mathcal{F}_0$, $\varepsilon>0$, there exists a
polynomial $p$ in $a_i$, $a_i^*$ such that $\Vert
b-p\Vert<\varepsilon$. We can decompose $p$ into sum of two
polynomials $p_1$, $p_2$ such that $p_1$ is a sum of those
homogeneous terms of $p$ where each $a_i$ appears so many times as
$a_i^*$ and $p_2=p-p_1$. Using the relations in the algebra we have
\[
p_1\in C^*(a_i^n (a_i^*)^n,\ (a_i^*)^n a_i^n,\ i=1,2,\
n\in\mathbb{Z}_{+}).
\]
As
\[
\int_0^{2\pi}\int_0^{2\pi} \rho_{\lambda_1,\lambda_2}(p_2)
d\lambda_1 d\lambda_2 =0,
\]
\begin{align*}
\Vert b-p_1\Vert &
=\frac{1}{(2\pi)^2}\Vert\int_0^{2\pi}\int_0^{2\pi}
\rho_{\lambda_1,\lambda_2}(b-p)
d\lambda_1 d\lambda_2\Vert\le\\
&\le \frac{1}{(2\pi)^2}\int_0^{2\pi}\int_0^{2\pi} \Vert\rho_{\lambda_1,\lambda_2}(b-p)\Vert d\lambda_1 d\lambda_2=\\
&= \frac{1}{(2\pi)^2}\int_0^{2\pi}\int_0^{2\pi} \Vert b-p\Vert
d\lambda_1 d\lambda_2=\Vert b-p\Vert<\varepsilon,
\end{align*}
and hence $b\in C^*(a_i^n (a_i^*)^n,\ (a_i^*)^n a_i^n,\ i=1,2,\
n\in\mathbb{Z}_{+})$. These arguments are standard in operator
algebras theory and we give them just for the completeness.

It follows from the relations between generators of
$\mathcal{E}(-1,0)$ and the description of $\mathcal{F}_0$ that
$\mathcal{F}_0$ is commutative. However, $\mathcal{F}_y$ is
non-commutative, since, for example, in  $\mathcal{E}(-1,y)$,
$a_1a_1^*a_2a_2^*\ne a_2a_2^*a_1a_1^*$ . Hence
$\mathcal{F}_0\not\simeq\mathcal{F}_y$ and therefore
$\mathcal{E}(-1,0)\not\simeq\mathcal{E}(-1,y)$ when $0<|y|<1$. To
complete the proof recall that $\mathcal{E}(-1,y)$ is isomorphic to
$M_2(\mathbb{C})$ if $|y|=1$ while $\mathcal{E}(-1,0)$ is not.
\end{proof}
Finally, let us consider the family of $C^*$-algebras,
$\mathcal{E}(\varepsilon)$, defined as enveloping for WCAR with
\[
a_2a_1+a_1a_2=y,\ \varepsilon\le\mid y\mid\le 1.
\]
It will be convenient for us to consider the polar decomposition of
$y$, $y=re^{i\phi}$. Using the results of Theorem \ref{theorem2} it
is easy to get the description of $\mathcal{E}(\varepsilon)$. Let
\[
D_{\varepsilon} = \{(r,x_1,x_2)\mid\ \varepsilon\le r\le 1,\ 0\le
x_1\le\frac{1-r}{2},\ 0\le x_2\le\frac{1}{2}\}.
\]
\begin{theorem}\label{theorem4}
For any\ $0<\varepsilon<1$ the $C^*$-algebra
$\mathcal{E}(\varepsilon)$ can be realized as fol-lows:
\begin{gather*}
\mathcal{E}(\varepsilon)=\bigl\{f\in C(D_{\varepsilon}\times \mathbf{T}^3\to
M_4(\mathbb{C}))\ |\\
f(r,0,x_2,e^{i\phi},e^{i\phi_1},e^{i\phi_2})=f(r,0,x_2,e^{i\phi},1,e^{i\phi_2}),
\\
f(r,x_1,0,e^{i\phi},e^{i\phi_1},e^{i\phi_2})=f(r,x_1,0,e^{i\phi},e^{i\phi_1},1), \\
f\bigl(r,(1-r)/2,x_2,e^{i\phi},e^{i\phi_1},e^{i\phi_2}\bigr)=
f\bigl(r,(1-r)/2,0,e^{i\phi},e^{i\phi_1},1\bigr)\in M_2(\mathbb{C})\otimes\mathbf{1},\\
f\bigl(1,0,x_2,e^{i\phi},e^{i\phi_1},e^{i\phi_2}\bigr)=f\bigl(1,0,0,e^{i\phi},1,1\bigr),\\
\nu^*_2(\phi_2)f(r,x_1,1/2,e^{i\phi},e^{i\phi_1},e^{i\phi_2})\nu_2(\phi_2)\in
M_2(\mathbb{C})\oplus M_2(\mathbb{C}), \forall\ \phi_2\in[0,2\pi]
\bigr\}.
\end{gather*}
\end{theorem}
\begin{proof}
The proof is evident.
\end{proof}
\begin{proposition}
For any $\varepsilon_1$, $\varepsilon_2$, $0<\varepsilon_i<1$,
$i=1,2$, $\mathcal{E}(\varepsilon_1)$ is isomorphic to
$\mathcal{E}(\varepsilon_2)$.
\end{proposition}
\begin{proof}
The required isomorphism is induced by the following homeomorphism
$\vartheta\colon D_{\varepsilon_1}\times\mathbf{T}^3\rightarrow
D_{\varepsilon_2}\times\mathbf{T}^3$:
\[
(r,x_1,x_2,e^{i\phi},e^{i\phi_1},e^{i\phi_2})\mapsto (\varepsilon_2+
\frac{1-\varepsilon_2}{1-\varepsilon_1}(r-\varepsilon_1),
\frac{1-\varepsilon_2}{1-\varepsilon_1}x_1,x_2,e^{i\phi},e^{i\phi_1},e^{i\phi_2}).
\]
\end{proof}
\section{Description of $\mathcal{E}(-1)$.}
Our goal in this section to describe the "global" $C^*$-algebra.
We put $y:= r e^{i\phi}$ and $r_1:=\frac {r}{1-2x_1}$. When $0\le
r\le 1$ and $0\le x_1\le \frac{1-r}{2}$ one has $0\le r_1\le 1$.
So, representation given by (\ref{cch}) has the form
\begin{align}\label{cchr1}
a_1 & =\left(%
\begin{array}{cc}
  0 & e^{i\phi_1}\sqrt{x_1} \\
  \sqrt{1-x_1} & 0\\
\end{array}%
\right)\otimes%
\left(%
\begin{array}{cc}
  1 & 0 \\
  0 & 1 \\
\end{array}%
\right),\ \phi_1\in [0,2\pi),\ \phi_2\in [0,2\pi)\nonumber\\  %
\  \nonumber \\
a_2 & =\sqrt{1-r_1^2}%
\left(%
\begin{array}{cc}
  1 & 0 \\
  0 & -1 \\
\end{array}%
\right)\otimes
\left(%
\begin{array}{cc}
  0 & e^{i\phi_2}\sqrt{x_2} \\
  \sqrt{1-x_2} & 0 \\
\end{array}%
\right)+  \\ \ \nonumber \\
&+r_1 e^{i\phi}\left(%
\begin{array}{cc}
  0 & \sqrt{1-x_1} \\
  -e^{-i\phi_1}\sqrt{x_1} & 0 \\
\end{array}%
\right)\otimes\left(%
\begin{array}{cc}
  1 & 0 \\
  0 & 1 \\
\end{array}%
\right);\nonumber
\end{align}
Consider
\[
D=\{(r_1,x_1,x_2)\ :\ 0\le x_i\le\frac{1}{2},\ i=1,2,\ 0\le r_1\le
1\}\times\mathbf{T}^3
\]
In this section we show that $\mathcal{E}(-1)$ is isomorphic to
the $C^*$-algebra of continuous matrix-functions on $D$ generated
by $a_1,\ a_2$ given by (\ref{cchr1}). To do so we need the
following auxiliary lemma.
\begin{lemma}\label{l12}
Let $\pi$ be a representation defined by (\ref{cchr1})
corresponding to a tuple
$(r_1,1/2,x_2,e^{i\phi},e^{i\phi_1},e^{i\phi_2})$, where $0<r_1\le
1$. Then $\pi$ is equivalent to representation $\widetilde{\pi}$
defined by (\ref{cchr1}) with parameters
$(0,1/2,t,1,e^{i\phi_1},e^{i\psi})$, where
\begin{align*}
& \sqrt{t(1-t)}=\mid z\mid,\ \psi=arg\ z,\\ &
z=e^{i\phi_2}(1-r_1^2)\sqrt{x_2(1-x_2)}-\frac{r_1^2
e^{i(2\phi-\phi_1)}}{2}
\end{align*}
\end{lemma}
\begin{proof}
First we note that the representation $\widetilde{\pi}$
corresponding to $(0,1/2,t,1,e^{i\phi_1},e^{i\psi})$ coincides
with a representation of $\mathcal{E}(0,-1)$ given by
(\ref{regg1}) with $x_1=1/2$, $\phi_1$ and $\phi_2=\psi$ and it is
equivalent to the direct sum of two-dimensional irreducible
representations given by formula (\ref{2a}), Section 3.2.:
\[
\widetilde{\pi}\sim \pi_1\oplus\pi_2,
\]
where $\pi_1$ corresponds to parameters $x_2=t$, $\phi_1/2$ and
$\phi_2=\psi$ and $\pi_2$ corresponds to
$(t,\frac{\phi_1}{2}+\pi,\psi)$. It is easy to see that
$\pi_1\not\sim\pi_2$ if $t\ne\frac{1}{2}$ and $\pi_1\sim\pi_2$
otherwise.

Thus, to prove lemma we have to show that for tuple
$(r_1,1/2,x_2,e^{i\phi},e^{i\phi_1},e^{i\phi_2})$, $0<r_1\le 1$
the corresponding representation $\pi$ is equivalent to
$\pi_1\oplus\pi_2$.

Since $a_1,\ a_2$ defined by $(\ref{cchr1})$ with
$x_1=\frac{1}{2}$ give a reducible representation of WCAR and
satisfy $a_1a_2+a_2a_1=0$, according to results of Sec. 3.2, $\pi$
is equivalent to the direct sum of two-dimensional irreducible
representations of $\mathcal{E}(-1,0)$, $\pi\sim
\widehat{\pi}_1\oplus\widehat{\pi}_2$. Further
\[
a_2^2=e^{i\phi_2}(1-r_1^2)\sqrt{x_2(1-x_2)}-\frac{r_1^2
e^{i(2\phi-\phi_1)}}{2}:=z
\]
implies that $\sigma(a_2^*a_2)=\{t,1-t\}$, where
$\sqrt{t(1-t)}=\mid z\mid$ and $u_2^2=e^{i\psi}$, $\psi=arg\ z$,
where $a_2=u_2(a_2^*a_2)^{\frac{1}{2}}$ is a polar decomposition.

If $\mid z\mid<\frac{1}{2}$ then $0\le t<\frac{1}{2}$, and taking
into account $a_1^2=\frac{e^{i\phi_1}}{2}$ we conclude that in
this case $\widehat{\pi}_1$, $\widehat{\pi}_2$ are determined by
formulas (\ref{2a}), Sec. 3.2, with $x_2=t$, $\phi_2=\psi$. To
describe the decomposition of $\pi$ completely it is remained to
verify whether or not $\widehat{\pi}_1$ and $\widehat{\pi}_2$ are
equivalent. To do so we compute the dimension of commutant of
$\pi$. It is a routine to verify that if $\mid z\mid<\frac{1}{2}$
then
\[
\dim \{\pi(a_i),\ i=1,2\}^{\prime}=2.
\]
Hence $\widehat{\pi}_1\not\sim\widehat{\pi}_2$ and we can suppose
that $\widehat{\pi}_1$, $\widehat{\pi}_2$ correspond to (\ref{2a})
with tuples $(t,\frac{\phi_1}{2},\psi)$ and
$(t,\frac{\phi_1}{2}+\pi,\psi)$ respectively. I.e., if $\mid
z\mid<\frac{1}{2}$, then $\widehat{\pi}_i\sim\pi_i$, $i=1,2$.

If $\mid z\mid=\frac{1}{2}$ then $t=\frac{1}{2}$,
$\psi=2\phi-\phi_1+\pi\mod 2\pi$, and
\begin{equation}\label{specq}
a_1^2=\frac{e^{i\phi_1}}{2},\
a_2^2=\frac{e^{i(2\phi-\phi_1+\pi)}}{2}.
\end{equation}
(Moreover if  $0<r_1<1$ we additionally have $x_2=\frac{1}{2}$ and
$\phi_2=2\phi-\phi_1+\pi\mod 2\pi$.)

Since the two-dimensional irreducible representation of
$\mathcal{E}(0,-1)$ with $a_1$, $a_2$ satisfying (\ref{specq}) is
unique we conclude that $\widehat{\pi}_1\sim\widehat{\pi}_2$ and
defined by (\ref{31}), Sec. 3.2, with $(\frac{\phi_1}{2},\psi)$.
Hence, when $\mid z\mid=\frac{1}{2}$ we also have
$\widehat{\pi}_i\sim\pi_i$, $i=1,2$  (recall, that the
representations given by (\ref{31}) corresponding to
$(\frac{\phi_1}{2},\psi)$ and $(\frac{\phi_1}{2}+\pi,\psi)$ are
equivalent).

So, we have shown that there exists a unitary matrix-function
\[
U=U(r_1,x_2,\phi,\phi_1,\phi_2),
\]
such that
\begin{align*}
U^*(r_1,x_2,\phi,\phi_1,\phi_2)&
a_i(r_1,1/2,x_2,e^{i\phi},e^{i\phi_1},e^{i\phi_2})
U(r_1,x_2,\phi,\phi_1,\phi_2)= \\ & =
a_i(0,1/2,t,1,e^{i\phi_1},e^{i\psi}),
\end{align*}
where $a_i$, $i=1,2$, are defined by (\ref{cchr1}) and $t$, $\psi$
are specified above.
\end{proof}
\begin{remark}
The matrix-function $U$ can be written explicitly, however we do
not give it here.
\end{remark}
Now we are ready to formulate the main results of this section.
\begin{theorem}
\[
\mathcal{E}(-1)\simeq
C^*(a_i(r_1,x_1,x_2,e^{i\phi},e^{i\phi_1},e^{i\phi_2}),\ i=1,2)
\]
where $a_i$, $i=1,2$, are given by (\ref{cchr1}) and
$(r_1,x_1,x_2,e^{i\phi},e^{i\phi_1},e^{i\phi_2})\in D$.
\end{theorem}
\begin{proof}
The result follows directly from the definition of enveloping
$C^*$-algebra and results of Theorems \ref{theorem2},
\ref{theorem3} and Lemma \ref{l12}.
\end{proof}
\begin{theorem} The $C^*$-algebra $\mathcal{E}(-1)$ is isomorphic
to the  $C^*$-algebra of continuous matrix-functions $f\colon
D\rightarrow M_4(\mathbb{C})$ satisfying the following boundary
conditions
\begin{gather*}
f(r_1,0,x_2,e^{i\phi},e^{i\phi_1},e^{i\phi_2})=
f(r_1,0,x_2,e^{i\phi},1,e^{i\phi_2})\\
f(r_1,x_1,0,e^{i\phi},e^{i\phi_1},e^{i\phi_2})=
f(r_1,x_1,0,e^{i\phi},e^{i\phi_1},1)\\
f(1,x_1,x_2,e^{i\phi},e^{i\phi_1},e^{i\phi_2})=
f(1,x_1,0,e^{i\phi},e^{i\phi_2},1)\in
M_2(\mathbb{C})\otimes\mathbf{1}_2\\
f(0,x_1,x_2,e^{i\phi},e^{i\phi_1},e^{i\phi_2})=
f(0,x_1,x_2,1,e^{i\phi_1},e^{i\phi_2})\\ \nu_2^*(\phi_2)
f(r_1,x_1,1/2,e^{i\phi},e^{i\phi_1},e^{i\phi_2})\nu_2(\phi_2)\in
M_2(\mathbb{C})\oplus M_2(\mathbb{C})\\ \nu_1^*(\phi_1,\phi_2)
f(0,1/2,x_2,e^{i\phi},e^{i\phi_1},e^{i\phi_2})\nu_1(\phi_1,\phi_2)\in
M_2(\mathbb{C})\oplus M_2(\mathbb{C})\\
U^*(r_1,x_2,\phi,\phi_1,\phi_2)
f(r_1,1/2,x_2,e^{i\phi},e^{i\phi_1},e^{i\phi_2})
U(r_1,x_2,\phi,\phi_1,\phi_2)=\\ =
f(0,1/2,t,1,e^{i\phi_1},e^{i\psi}),
\end{gather*}
where $t$, $\psi$ and $U(\cdot)$ are specified in Lemma \ref{l12}.
\end{theorem}
\begin{proof}
The proof is the same as for Theorems \ref{theorem2},
\ref{theorem3}.
\end{proof}
\begin{remark}
Note that
$U(0,x_2,\phi,\phi_1,\phi_2)=U(0,x_2,1,\phi_1,\phi_2)=\mathbf{1}_4$.
\end{remark}
\noindent {\bf Acknowledgements.} We are indebted to Prof. P.E.T.
J{\o}rgensen and Prof. Yu.S. Samo\u\i{}lenko for helpful discussions
on the subject.
\bibliographystyle{amsalpha}

\end{document}